\documentclass{amsart}
\usepackage{amssymb,amsfonts,latexsym}

\newtheorem{theorem}{Theorem}[section]

\newtheorem{proposition}[theorem]{Proposition}
\newtheorem{corollary}[theorem]{Corollary}
\theoremstyle{definition}

\newtheorem{remark}[theorem]{Remark}

\newcommand{\ben}{\begin{enumerate}}
\newcommand{\een}{\end{enumerate}}

\newcommand{\tensor}{\otimes}

\newcommand{\directsum}{\oplus}
\newcommand{\Epsilon}{\mathcal{E}}

\begin{document}

\title[On the Hopf-Schur group of a field]
{On the Hopf-Schur group of a field}

\author{Eli Aljadeff}
\address{Department of Mathematics, Technion-Israel Institute of
Technology, Haifa 32000, Israel}
\email{aljadeff@tx.technion.ac.il}
\author{Juan Cuadra}
\address{Departamento de \'Algebra y An\'alisis Matem\'atico, Universidad de Almer\'{\i}a, E04120 Almer\'{\i}a, Spain}
\email{jcdiaz@ual.es}
\author{Shlomo Gelaki}
\address{Department of Mathematics, Technion-Israel Institute of
Technology, Haifa 32000, Israel}
\email{gelaki@math.technion.ac.il}
\author{Ehud Meir}
\address{Department of Mathematics, Technion-Israel Institute of
Technology, Haifa 32000, Israel}
\email{ehudm@techunix.technion.ac.il}

\date{August 14, 2007}

\begin{abstract}
Let $k$ be any field. We consider the Hopf-Schur group of $k$,
defined as the subgroup of the Brauer group of $k$ consisting of
classes that may be represented by homomorphic images of Hopf
algebras over $k$. We show here that twisted group algebras and
abelian extensions of $k$ are quotients of cocommutative and
commutative Hopf algebras over $k$, respectively. As a consequence
we prove that any tensor product of cyclic algebras over $k$ is a
quotient of a Hopf algebra over $k$, revealing so that the
Hopf-Schur group can be much larger than the Schur group of $k$.
\end{abstract}

\maketitle

\begin{section}{Introduction}

Let $k$ be any field and let $H$ be a finite dimensional Hopf
algebra over $k$. Let $C$ be a simple algebra which is a
homomorphic image of $H$. Clearly, $C$ determines an element in
$Br(L)$, the Brauer group of $L$, where $L$ is the center of $C$.
We say that an $L$-central simple algebra $C$ is {\em Hopf-Schur}
over $L$ if it is a homomorphic image of a finite dimensional Hopf
algebra over $k$, where $k$ is a subfield of $L$. The Hopf-Schur
group $HS(L)$ is the subgroup of $Br(L)$ generated by (and in fact
consisting of) classes that may be represented by Hopf-Schur
algebras over $L$. Since any Hopf-Schur algebra over $L$ is a
homomorphic image of a Hopf algebra over $L$ (just by extension of
scalars) we may restrict our discussion to Hopf-Schur algebras
over $k$ which are homomorphic images of finite dimensional Hopf
algebras over $k$. Observe that a $k$-central simple algebra is a
homomorphic image of a {\em semisimple} $k$-Hopf algebra $H$ if
and only if it is a simple component of the Wedderburn
decomposition of $H$. The group $HS(k)$ clearly contains $S(k)$,
the Schur group of $k$. Recall that a $k-$central simple algebra
is said to be {\em Schur} over $k$ if it is a homomorphic image of
a group algebra $kG$ for some finite group $G$. The group $S(k)$
is the subgroup of $Br(k)$ consisting of (Brauer) classes that may
be represented by Schur algebras over $k$. In general $S(k)$ is
``very small" compared to the full Brauer group. For instance if
$k$ contains no non-trivial cyclotomic extensions then $S(k)=0$
whereas $Br(k)$ may be large (e.g.,
$k=\mathbb{C}(x_{1},..,x_{n}),\; n\geq 2$, see [FS]). Although the
class of semisimple Hopf algebras over $k$ is much richer than the
class of group algebras over $k$, no example was known of an
element in $HS(k)$ which is not in $S(k)$ (we refer the reader to
[J] and [Y] for a comprehensive account on the Schur group). The
main goal of this article is to show that such examples are
abundant. In fact we will show that any cyclic algebra (and hence
any tensor product of cyclic algebras) is a homomorphic image of a
Hopf algebra. Recall that a {\em cyclic algebra} $C$ over $k$ is a
crossed product algebra $(L/k, G, \alpha)$, $\alpha\in
Z^2(G,L^{*})$, where $G:=Gal(L/k)$ is a cyclic group and $L^*$
denotes the multiplicative group of $L$. The main result of the
paper is the following:

\begin{theorem}\label{main1}
Let $C$ be a $k-$central simple algebra which is isomorphic to a
tensor product of cyclic algebras over $k$. Then $C$ is a
homomorphic image of a semisimple Hopf algebra over $k$.
\end{theorem}

Recall that it is a major open problem in Brauer group theory
whether any element in $Br(k)$ is Brauer equivalent to a product
of cyclic algebras. This is known to be true for number fields
(Albert-Brauer-Hasse-Noether Theorem) and for fields containing
all roots of unity (Merkurjev-Suslin Theorem). Thus, we naturally
raise the question whether $HS(k)=Br(k)$ for arbitrary fields. Of
course, this question may be solved (in the positive)
independently of the above mentioned problem on cyclic algebras.

The first step in the proof of Theorem \ref{main1} is carried out
in Section 2 where we prove the following result, which is of
independent interest.

\begin{theorem}\label{main2}
Let $G$ be a finite group, and let $[\alpha]\in H^2(G,k^*)$ be an
element of order $m$.
\begin{enumerate}
\item
The twisted group algebra $k^{\alpha}G$ is a homomorphic image of
a finite dimensional cocommutative Hopf algebra $A$ over $k$ of
dimension $m|G|$. When the characteristic of $k$ does not divide
$|G|$, $A$ is semisimple.
\item
Suppose the characteristic of $k$ does not divide $m$. Then the
$k-$Hopf algebra $A$ is a form of the group algebra
$\bar{k}\widehat{G}$, i.e., $A\otimes_k \bar{k}\cong
\bar{k}\widehat{G}$, where $\widehat{G}$ is the central extension
(induced by $[\alpha]$) of $G$ by $\mathbb{Z}/m\mathbb{Z}$.
\end{enumerate}
\end{theorem}

\begin{remark}
We shall explain how $[\alpha]$ determines an extension of $G$ by
$\mathbb{Z}/m\mathbb{Z}$ in Section \ref{twisted-algebras}.
\end{remark}

Now let $C=(L/k, G, \alpha)$ be a cyclic algebra. It is well known
that by passing to a cohomologous $2$-cocycle we can assume that
$\alpha$ has values in $k^{*}$ (rather than in $L^{*}$). Hence $C$
contains a $k$-subalgebra isomorphic to the twisted group algebra
$k^{\alpha}G$. The next result, accomplished in Section 3, is the
second step in the proof of Theorem \ref{main1}.

\begin{theorem}\label{main3}
Let $L/k$ be a finite abelian extension with $G:=Gal(L/k)$ and let
$\textrm{Fun}(\mathbb{Z}/2\mathbb{Z}\ltimes G,k)$ be the $k-$Hopf
algebra of functions on $\mathbb{Z}/2\mathbb{Z}\ltimes G$ with
values in $k$, where $\mathbb{Z}/2\mathbb{Z}$ acts on $G$ by
inversion. Then there exists a commutative semisimple Hopf algebra
$H$ over $k$ of dimension $2[L:k]$ which satisfies the following:
\begin{enumerate}
\item
$L$ is a homomorphic image of $H$.
\item
$H\otimes_k L\cong \textrm{Fun}(\mathbb{Z}/2\mathbb{Z}\ltimes G,k)
\otimes_k L \cong \textrm{Fun}(\mathbb{Z}/2\mathbb{Z}\ltimes
G,L)$. In particular, $H$ is a $k-$form of
$\textrm{Fun}(\mathbb{Z}/2\mathbb{Z}\ltimes G,k)$.

\end{enumerate}
\end{theorem}

Thus we have constructed two Hopf algebras $A$ and $H$ over $k$,
such that $k^{\alpha}G$ and $L$ are quotients of them,
respectively. In order to complete the proof of Theorem
\ref{main1} we ``amalgamate'' $A$ and $H$ to obtain a Hopf algebra
$X$ over $k$ of which $C=(L/k, G, \alpha)$ is a quotient. This is
done in Section \ref{sec-main-thm}.


In 1978 Lorenz and Opolka introduced the {\em projective Schur}
subgroup, denoted by $PS(k)$, of $Br(k)$. It consists of Brauer
classes that may be represented by homomorphic images of twisted
group algebras, i.e., algebras of the form $k^{\alpha}G$ where
$\alpha$ is a $2$-cocycle of $G$ with coefficients in $k^{*}$.
Clearly $PS(k)$ contains $S(k)$. Moreover, it is easy to construct
examples (e.g. symbol algebras $(a,b)_{n}$) which are in $PS(k)$
but not in $S(k)$. In general a twisted group algebra
$k^{\alpha}G$ is not a Hopf algebra over $k$. In fact it is not
difficult to show that $k^{\alpha}G$ is a Hopf algebra over $k$ if
and only if $\alpha$ is cohomologically trivial. As an immediate
consequence of Theorem \ref{main2} we have:

\begin{corollary}
$HS(k)\supseteq PS(k)$.
\end{corollary}

\begin{remark}
We point out that there are fields $k$ for which $HS(k)$ properly
contains $PS(k)$. This follows from Theorem \ref{main1} and the
fact that there are cyclic algebras over suitable fields $k$ (e.g.
$k=\mathbb{Q}(x)$) which do not represent elements in $PS(k)$ (see
[AS2]). On the other hand, there exist fields $k$ for which
$HS(k)=PS(k)$. For example, for number fields $PS(k)=Br(k)$, so
$PS(k)=HS(k)=Br(k).$
\end{remark}

\begin{remark} We do not know of any $k$-central simple
algebra which is not a homomorphic image of a semisimple Hopf
algebra. Also, from the above, it is not clear whether any element
in $HS(k)$ is split by an abelian extension of $k$. The same
question for $PS(k)$ was answered in the affirmative in [AS1].
\end{remark}

\end{section}

\begin{section}{The proof of Theorem \ref{main2}}\label{twisted-algebras}

Let $k$ be any field, and let $G$ be a finite group. Recall that
the twisted group algebra $k^{\alpha}G$, with respect to a
2-cocycle $\alpha\in Z^2(G,k^*)$, is the $k$-algebra spanned as a
$k$-vector space by the elements $U_{\sigma}$, $\sigma\in G$, with
multiplication defined by the formula
$$U_{\sigma}U_{\tau}=\alpha(\sigma,\tau)U_{\sigma\tau},\;\sigma,\tau\in G.$$
It is well known that up to isomorphism of $k$-algebras,
$k^{\alpha}G$ does not depend on the cocycle $\alpha$ but only on
the cohomology class $[\alpha] \in H^2(G,k^*)$ it represents. In
this section we exhibit a Hopf algebra $A$ over $k$ such that
$k^{\alpha}G$ is a quotient of $A$. We will first define the
algebra structure on $A$ and then we will define the counit,
comultiplication and antipode.

Recall that every element in the group $H^2(G,k^*)$ is of finite
order. Therefore there exist a natural number $m$ and a 1-cochain
$f$ such that $\alpha^m=df$. This means that
\begin{equation}
\label{alpha and
f}\alpha^m(\sigma,\tau)=f(\sigma)f(\tau)/f(\sigma\tau),
\;\sigma,\tau\in G.
\end{equation}
Without loss of generality we may assume that
$\alpha(\sigma,1)=\alpha(1,\sigma)=f(1)=1$ for every $\sigma \in
G$. From now on we shall denote by $m$ the order of $[\alpha]$
For $n=0,\ldots,m-1$, we shall denote the basis of $k^{\alpha^n}G$
by $\{U^{(n)}_{\sigma}\}_{\sigma\in G}$. Consider the algebra
$$A:=\bigoplus_{n=0}^{m-1}k^{\alpha^n}G.$$ The algebra $A$ has
a counit $\epsilon$ given by
$\epsilon(U^{(n)}_{\sigma})=\delta_{n,0}$. It is easy to see that
$\epsilon$ is an algebra morphism. Notice that the algebra $A$
will not be semisimple if the characteristic of $k$ divides the
order of $G$ (since $kG$ is a direct summand of $A$). On the other
hand, if $char(k)$ does not divide $\vert G \vert$, then
$k^{\alpha^i}G$ is semisimple for all $i=0,...,m-1$. This shows
that $A$ is semisimple if and only if $char(k)$ does not divide
$\vert G \vert$. We shall now make $A$ into a Hopf algebra. For
every $0\leq r,l< m$, $\sigma\in G$, define
$$\xi^{\sigma}_{r,l}=\left\{\begin{array}{ll}
1 & \textrm{if}\; r+l< m \\
\frac{1}{f(\sigma)} & \textrm{if}\; r+l\geq m
\end{array}\right.$$ Note that
$\xi^{\sigma}_{r,l}=\xi^{\sigma}_{l,r}$ for all $r$, $l$ and
$\sigma$.

Now define a map $\Delta :A\to A\otimes A$ by
$$\Delta(U^{(n)}_{\sigma})=\sum_{{r+l\equiv n}\atop {mod. m}}{\xi^{\sigma}_{r,l}U^{(r)}_{\sigma}\tensor U^{(l)}_{\sigma}},$$
and a map $S:A\to A$ by
$$S(U^{(n)}_{\sigma})=\left\{\begin{array}{ll}
U^{(0)}_{\sigma^{-1}} & \textrm{if}\; n=0 \vspace{3pt} \\
\frac{\alpha^n(\sigma,\sigma^{-1})}{f(\sigma^{-1})}U^{(m-n)}_{\sigma^{-1}}
& \textrm{if}\; n > 0
\end{array}\right.$$
Note that $\Delta$ is invariant under the flip map $A\otimes A \to
A\otimes A$. Therefore, after we will prove that $A$ is a Hopf
algebra, it will follow that $A$ is cocommutative.

\begin{proposition}
The above formulas define a cocommutative Hopf algebra structure
on the algebra $A$.
\end{proposition}

\begin{proof}
We shall show, by a direct calculation, that all the axioms of a
Hopf algebra are valid in $A$.

{\bf 1.} $\Delta$ is coassociative. Since

$$(\Delta\tensor 1)\Delta(U^{(n)}_{\sigma}) = (\Delta\tensor
1)(\hspace{-4pt} \sum_{r+l\equiv
n}{\xi^{\sigma}_{r,l}U^{(r)}_{\sigma}\tensor U^{(l)}_{\sigma}})=
\hspace{-7pt} \sum_{z+q+l\equiv
n}{\xi^{\sigma}_{z,q}\xi^{\sigma}_{z+q,l}U^{(z)}_{\sigma}\tensor
U^{(q)}_{\sigma}\tensor U^{(l)}_{\sigma}}$$ and $$(1\tensor
\Delta)\Delta(U^{(n)}_{\sigma}) = (1\tensor \Delta)(\hspace{-4pt}
\sum_{r+l\equiv n}{\xi^{\sigma}_{r,l}U^{(r)}_{\sigma}\tensor
U^{(l)}_{\sigma}}) =\hspace{-7pt} \sum_{r+s+t\equiv
n}{\xi^{\sigma}_{r,s+t}\xi^{\sigma}_{s,t}U^{(r)}_{\sigma}\tensor
U^{(s)}_{\sigma}\tensor U^{(t)}_{\sigma}},$$ we need to prove that
for every $r,s,t<m$ we have
\begin{equation}\label{xi-condition1}\xi^{\sigma}_{r,s+t}\xi^{\sigma}_{s,t} =
\xi^{\sigma}_{r,s}\xi^{\sigma}_{r+s,t}\end{equation} where the sum
of the indices is taken modulo $m$. We prove this by considering
the following three possible cases:

\begin{itemize}
\item[a.] $r+s+t<m$: In this case both sides of the equation equal 1.

\item[b.] $m\leq r+s+t<2m$: In this case in both sides of the equation
one of the terms will equal $1/f(\sigma)$ and the other will equal
1.

\item[c.] $2m\leq r+s+t$: In this case all the terms in the equation
will be equal to $1/f(\sigma)$ and so both sides of the equation
will equal $1/f(\sigma)^2$.
\end{itemize}

{\bf 2.} $\Delta$ is an algebra map. First,
$$\Delta(1_A)=\Delta(\sum_{n=0}^{m-1}U^{(n)}_{1})=
\sum_{r,l=0}^{m-1}{\xi^1_{r,l}U^{(r)}_1\tensor U^{(l)}_1} =
1_{A\tensor A}.$$ So $\Delta$ sends the unit of $A$ to the unit of
$A\tensor A$. Next, if $n \neq p$, then
$$\begin{array}{ll}
\Delta(U^{(n)}_{\sigma})\Delta(U^{(p)}_{\tau}) & = ({\displaystyle
\sum_{r+l\equiv n}}{\xi^{\sigma}_{r,l}U^{(r)}_{\sigma}\tensor
U^{(l)}_{\sigma}})({\displaystyle \sum_{a+b\equiv
p}}{\xi^{\tau}_{a,b}U^{(a)}_{\tau}\tensor U^{(b)}_{\tau}}) \\
& = 0 \\
& =\Delta(0) \\
& =\Delta(U^{(n)}_{\sigma}U^{(p)}_{\tau}),
\end{array}$$
and
$$\begin{array}{ll}
\Delta(U^{(n)}_{\sigma})\Delta(U^{(n)}_{\tau}) & = ({\displaystyle
\sum_{r+l\equiv n}}{\xi^{\sigma}_{r,l}U^{(r)}_{\sigma}\tensor
U^{(l)}_{\sigma}})(\sum_{r+l\equiv
n}{\xi^{\tau}_{r,l}U^{(r)}_{\tau}\tensor U^{(l)}_{\tau}}) \vspace{3pt} \\
& = {\displaystyle \sum_{r+l\equiv
n}}{\xi^{\sigma}_{r,l}\xi^{\tau}_{r,l}\alpha^r(\sigma,\tau)U^{(r)}_{\sigma\tau}\tensor
\alpha^l(\sigma,\tau)U^{(l)}_{\sigma\tau}} \vspace{3pt} \\
& = {\displaystyle \sum_{r+l\equiv
n}}{\xi^{\sigma}_{r,l}\xi^{\tau}_{r,l}\alpha^{r+l}(\sigma,\tau)U^{(r)}_{\sigma\tau}\tensor
U^{(l)}_{\sigma\tau}}.
\end{array}$$
On the other hand,
$$\Delta(U^{(n)}_{\sigma}U^{(n)}_{\tau})=\Delta(\alpha^n(\sigma,\tau)U^{(n)}_{\sigma\tau})=
\sum_{r+l\equiv
n}{\xi^{\sigma\tau}_{r,l}\alpha^n(\sigma,\tau)U^{(r)}_{\sigma\tau}\tensor
U^{(l)}_{\sigma\tau}}.$$ So we need to prove that if $r+l\equiv n$
then
\begin{equation}\label{xi-condition2}\xi^{\sigma}_{r,l}\xi^{\tau}_{r,l}\alpha^{r+l}(\sigma,\tau)=
\xi^{\sigma\tau}_{r,l}\alpha^n(\sigma,\tau) .\end{equation} There
are two cases here:

\begin{itemize}
\item[a.] $r+l=n$: In this case both sides of the equation equal
$\alpha^n(\sigma,\tau)$.

\item[b.] $r+l=m+n$: In this case the equation becomes
$1/(f(\sigma)f(\tau))\alpha^{n+m}(\sigma,\tau)=1/f(\sigma,\tau)\alpha^n(\sigma,\tau)$
which is equivalent to
$\alpha^m(\sigma,\tau)=f(\sigma)f(\tau)/f(\sigma\tau)$, which is
exactly equation \ref{alpha and f}.
\end{itemize}

{\bf 3.}  $\Delta$ and $\epsilon$ make $A$ into a coalgebra. We
compute
$$(\epsilon\tensor 1)\Delta(U^{(n)}_{\sigma})=(\epsilon\tensor
1)(\sum_{r+l\equiv n}{\xi^{\sigma}_{r,l}U^{(r)}_{\sigma}\tensor
U^{(l)}_{\sigma}})=\xi^{\sigma}_{0,n}U^{(n)}_{\sigma}=
U^{(n)}_{\sigma}$$ and the result follows. The fact that
$(1\tensor \epsilon)\Delta=id_A$ follows by a similar argument or
by the above argument and the fact that $A$ is cocommutative.

{\bf 4.} $S$ is an antipode for $A$. If $0<n<m$ then
$$\begin{array}{ll}
m(S\tensor 1)\Delta(U^{(n)}_{\sigma}) & = m(S\tensor
1)({\displaystyle \sum_{r+l\equiv n}}
{\xi^{\sigma}_{r,l}U^{(r)}_{\sigma}\tensor U^{(l)}_{\sigma}}) \vspace{3pt} \\
& = {\displaystyle \sum_{r+l\equiv
n}}{\xi^{\sigma}_{r,l}\alpha^r(\sigma,\sigma^{-1})/f(\sigma^{-1})U^{(m-r)}_{\sigma^{-1}}U^{(l)}_{\sigma}} \vspace{3pt} \\
& = 0  \vspace{3pt} \\
& =\epsilon(U^{(n)}_{\sigma}).
\end{array}$$
The third equality follows from the fact that if $m-r=l$, then
$m=r+l$, contradicting the assumption that $r+l\equiv n$. If $n=0$
then
$$\begin{array}{ll}
m(S\tensor 1)\Delta(U^{(0)}_{\sigma}) & \hspace{-5pt} = m(S\tensor
1)\big(U^{(0)}_{\sigma}\tensor U^{(0)}_{\sigma} + {\displaystyle
\sum_{r=1}^{m-1}}{\xi^{\sigma}_{r,m-r}U^{(r)}_{\sigma}\tensor
U^{(m-r)}_{\sigma}}\big) \vspace{3pt} \\
& \hspace{-5pt} =
U^{(0)}_{\sigma}U^{(0)}_{\sigma^{-1}}+{\displaystyle
\sum_{r=1}^{m-1}}{\xi^{\sigma}_{r,m-r}
\alpha^r(\sigma,\sigma^{-1})/f(\sigma^{-1}) U^{(m-r)}_{\sigma^{-1}}U^{(m-r)}_{\sigma}}) \vspace{3pt} \\
& \hspace{-5pt} = U^{(0)}_{1} \hspace{-2.5pt}+ \hspace{-2.5pt}
{\displaystyle \sum_{r=1}^{m-1}}{(1/f(\sigma))
(\alpha^r(\sigma,\sigma^{-1})/f(\sigma^{-1}))\alpha^{m-r}(\sigma^{-1},\sigma)
U^{(m-r)}_{1}} \vspace{3pt} \\
& \hspace{-5pt} =
U^{(0)}_1+\alpha^m(\sigma,\sigma^{-1})/(f(\sigma)f(\sigma^{-1}))
{\displaystyle \sum_{r=1}^{m-1}}{U^{(m-r)}_1} \vspace{3pt} \\
& \hspace{-5pt} = {\displaystyle \sum_{r=0}^{m-1}} {U^{(r)}_1} \vspace{3pt} \\
& \hspace{-5pt} =1_A \vspace{3pt} \\
& \hspace{-5pt} =\epsilon(U^{(0)}_{\sigma}),
\end{array}$$
and therefore $S$ is an antipode. During the computation we have
used in the third equality the fact that
$\alpha(\sigma,\sigma^{-1})=\alpha(\sigma^{-1},\sigma)$, and in
the fourth one that
$$\alpha^m(\sigma,\sigma^{-1})=f(\sigma)f(\sigma^{-1})/f(1)=f(\sigma)f(\sigma^{-1}).$$
This completes the proof that $A$ is a Hopf algebra. \qed

\begin{remark}
It is well known that for representations $V$ and $W$ of
$k^{\alpha^i}G$ and $k^{\alpha^j}G$, respectively, $V\otimes W$ is
a representation of the twisted group algebra $k^{\alpha^{i+j}}G$.
This defines a tensor structure on $\textrm{Rep}(A)$, which
explains the comultiplication on $A$. It is also known that $V^*$
is a representation of $k^{\alpha^{-i}}G\cong k^{\alpha^{m-i}}G$.
This means that the category $\textrm{Rep}(A)$ is rigid, which
explains the antipode on $A$. The trivial one dimensional
representation of $A$, $k$, on which $U^{(n)}_{\sigma}$ acts as
the scalar $\delta_{n,0}$, explains the counit of $A$. Thus the
maps which make the algebra $A$ into a Hopf algebra can be
reconstructed from the rigid tensor structure of the category
$\textrm{Rep}(A)$. (See e.g., \cite{G} for more details on
reconstruction.)
\end{remark}

The twisted group algebra $k^{\alpha}G$ and every quotient of it
is clearly a quotient of $A$. Thus, we have showed that every
twisted group algebra, and hence every projective Schur algebra,
is a quotient of a Hopf algebra.

We shall now prove the second part of Theorem \ref{main2}.
Consider the cocycle $\alpha$ as a cocycle with values in
$\bar{k}^*$ rather than in $k^*$. Since $\alpha^m$ is
cohomologically trivial, there is a 1-cochain $f:G\rightarrow
\bar{k}^*$ such that
$$\alpha^m(\sigma,\tau)=f(\sigma)f(\tau)/f(\sigma\tau).$$
Define a 1-cochain $g:G\rightarrow \bar{k}^*$ by
$g(\sigma)=f(\sigma)^{1/m}$, where $f(\sigma)^{1/m}$ is any
element of $\bar{k}^*$ whose $m$th power equals $f(\sigma)$. It is
easy to see that all the values of the cocycle
$\bar{\alpha}=\alpha/dg$ are $m$th roots of unity, and that
$\bar{\alpha}$ is cohomologous to $\alpha$. Suppose that the
characteristic of $k$ does not divide $m$. In this case the group
of $m$th roots of unity in $\bar{k}$ is isomorphic to
$\mathbb{Z}/m\mathbb{Z}$, and we can consider $\bar{\alpha}$ as a
central extension
$$\bar{\alpha} : 1\rightarrow \mathbb{Z}/m\mathbb{Z}
\rightarrow \widehat{G}\rightarrow G \rightarrow 1.$$

Let $\zeta\in \bar{k}$ be a primitive $m$th root of unity. We have
seen that there is a cocycle $\bar{\alpha}$ cohomologous to
$\alpha$ such that $\bar{\alpha}(\sigma,\tau)\in
\langle\zeta\rangle\cong \mathbb{Z}/m\mathbb{Z}$ for every
$\sigma,\tau\in G$. So we shall assume henceforth that $\alpha$
has all its values in $\langle\zeta\rangle$. In that case
$\alpha^m$ is the cocycle which is identically $1$. The map $f$ is
then taken as that one being identically $1$. The group
$\widehat{G}$ is the set of all ordered pairs of the form
$(\sigma,i)$, where $\sigma\in G$ and $i\in
\mathbb{Z}/m\mathbb{Z}$, with multiplication given by the formula
$$(\sigma,i)(\tau,j)=(\sigma\tau,i+j+\alpha(\sigma,\tau)).$$ We
consider $\alpha$ as a cocycle with values in the group
$\mathbb{Z}/m\mathbb{Z}$, and define the following
$\bar{k}$-linear map $$\phi:\bar{k}\widehat{G}\rightarrow
\bar{k}\otimes_k A,\;\;(\sigma,i)\mapsto
\sum_{j=0}^{m-1}{\zeta^{ij}\otimes U^{(j)}_{\sigma}}.$$ A
straightforward verification shows that this is a Hopf algebra
isomorphism, and thus $A$ is a $\bar{k}-$form of
$\bar{k}\widehat{G}$, as desired. This completes the proof of
Theorem \ref{main2}.
\end{proof}

\end{section}

\begin{section}{The proof of Theorem \ref{main3}}\label{section-H}
Let $L$ be an abelian Galois extension of $k$ with (abelian)
Galois group $G$ of order $n$. Recall the definition of
$\textrm{Fun}(G,k)$, the algebra of functions on $G$ with values
in $k$. It is the dual of the group algebra $kG$. It has a
$k$-basis consisting of the mutually orthogonal idempotents
$\{e_{\sigma}\}_{\sigma \in G}$ given by
$e_{\sigma}(\tau)=\delta_{\sigma,\tau}$. We define $H$, as an
algebra, to be
$$H:=L\directsum \textrm{Fun}(G,k).$$ We shall denote the unit of $L$ in $H$ by $1_L$,
to avoid confusion with the unit of $H$. Let us write the unit of
$Fun(G,k)$ as $e=\sum_{\sigma \in G}{e_{\sigma}}$. The unit of $H$
will then be $1_L+e$. There is also a natural counit $\epsilon$,
given by $\epsilon(e_{\sigma})=\delta_{\sigma, 1}$, and
$\epsilon(L)=0$. It is easy to see that $\epsilon$ is an algebra
map. $H$ also has an anti-algebra morphism, $S:H\rightarrow H$,
defined by $S(x)=x$ for $x\in L$, and
$S(e_{\sigma})=e_{\sigma^{-1}}$. This will be the antipode of $H$.
Notice that $H$ is a semisimple commutative algebra, and that $L$
is a quotient of $H$. We will now define the comultiplication on
$H$.

Since $L$ is a direct summand of $H$, $L\tensor _k L$ will be a
direct summand of $H\tensor _k H$, so we begin with some analysis
of $L\tensor _k L$. From Galois theory we know that $L\tensor _k
L$ decomposes as the direct sum of $n$ copies of $L$. Explicitly,
for every $\sigma\in G$ consider the algebra map
$\Psi_{\sigma}:L\tensor _k L \rightarrow L$ defined by
$\Psi_{\sigma}(x\tensor y)=\sigma(x)y$. Then we can write
$$L\tensor _k L = \directsum_{\sigma\in G}L\Epsilon_{\sigma}$$
where $\Epsilon_{\sigma}$ is the idempotent lying in the kernels
of all the maps $\Psi_{\tau}, \tau\neq\sigma$. For every $\sigma
\in G$, let $\omega_{\sigma}$ be the automorphism of $L\tensor L$
given by $\omega_{\sigma}(x\tensor y)=\sigma(x)\tensor\sigma(y)$.
Then for $\mu\in G$ we have
$\omega_{\sigma}(\Epsilon_{\mu})=\Epsilon_{\sigma\mu\sigma^{-1}}=
\Epsilon_{\mu}$(because $G$ is abelian). If we write
$\Epsilon_{\mu}=\sum_i{x_i\tensor y_i}$ then this means that
$\sum_i{x_i\tensor y_i}=\sum_i{\sigma(x_i)\tensor \sigma(y_i)}$.

Now define a map $\Delta :H\to H\otimes H$ by
$$\Delta(x)=\sum_{\mu\in G}{\big(\mu^{-1}(x)\tensor e_{\mu} + e_{\mu}\tensor \mu(x)\big)},\;x\in L,$$
and
$$\Delta(e_{\xi})=\sum_{\sigma\mu=\xi}{e_{\sigma}\tensor e_{\mu}} + \Epsilon_{\xi}.$$
\begin{proposition}\label{h}
The maps $\Delta$, $\epsilon$ and $S$ defined above equip $H$ with
the structure of a semisimple commutative Hopf algebra.
\end{proposition}

\begin{proof}
One can verify easily that $\Delta$ is an algebra map. The fact
that $\Delta$ is coassociative can be proved directly or by the
following argument. Consider the set $D:=Hom_{k-alg}(H,L)$; it
contains the following elements: for every $\sigma\in G$ we have
the map $\phi_{\sigma}:H\rightarrow L$ given by
$\phi_{\sigma}(x)=0$ for $x\in L$ and $\phi_{\sigma}(e_{\tau}) =
\delta_{\sigma,\tau}$, and the map $\zeta_{\sigma}:H\rightarrow L$
given by $\zeta_{\sigma}(x)=\sigma(x)$ for $x\in L$, and
$\zeta_{\sigma}(e_{\tau})=0$. It is easy to check that these $2n$
maps are all the $k$-algebra morphisms from $H$ to $L$. Since $L$
is commutative and $\Delta$ is multiplicative, the map $\Delta$
defines a convolution product * on $D$. It is easy to check that
the following holds by using the formulas for $\Delta$ given
above:
$$\phi_{\sigma}*\phi_{\tau}=\phi_{\sigma\tau}$$
$$\phi_{\sigma}*\zeta_{\tau}=\zeta_{\sigma\tau}$$
$$\zeta_{\tau}*\phi_{\sigma}=\zeta_{\tau\sigma^{-1}}$$
$$\zeta_{\sigma}*\zeta_{\tau}=\phi_{\tau^{-1}\sigma}.$$

We have used in the last equality that $$m_L(\zeta_{\sigma}\otimes
\zeta_{\tau})(\Epsilon_{\mu})=\Psi_{\sigma\tau^{-1}}(\tau\otimes
\tau)(\Epsilon_{\mu})=\Psi_{\sigma\tau^{-1}}(\Epsilon_{\mu}).$$ By
a direct verification, $($D$,*)$ is isomorphic as a set with a
binary operation to the group $\mathbb{Z}_2\ltimes G$, where
$\mathbb{Z}_2$ acts on $G$ by inversion. The correspondence
between these two sets is given by $\phi_{\sigma}\leftrightarrow
(0,\sigma)$, $\zeta_{\sigma}\leftrightarrow (1,\sigma)$, where we
write $\mathbb{Z}_2=\{0,1\}$. In particular, this means that $*$
is \textit{associative} on $D$, and therefore, for every 3
elements $x,y,z\in D$ we have
\begin{equation}\label{associativity} m_L(1\tensor m_L)(x\tensor
y\tensor z)(1\tensor \Delta)\Delta = m_L(m_L\tensor 1)(x\tensor
y\tensor z)(\Delta\tensor 1)\Delta.\end{equation} It follows that
Equation \ref{associativity} will hold if we replace $x,y,z$ by
$L$-linear combinations of elements of $D$. From Galois theory we
know that the elements of $G$ are $L$-linear independent. It
follows that the elements of $D$ are $L$-linear independent when
we consider them as elements of the $L$-vector space
$V=Hom_k(H,L)$, where $L$ acts by $(xf)(h)=xf(h)$, $x\in L$, $f\in
V$, and $h\in H$. Since $dim_LV=2n$, $D$ is a basis of $V$, and
therefore for every $x,y,z\in V$, Equation \ref{associativity}
holds. But this implies that $\Delta$ is coassociative. The fact
that $\epsilon$ is a counit follows from the following
computations ($x\in L$):
$$(\epsilon\tensor 1)\Delta(x)= (\epsilon\tensor 1)(\sum_{\mu\in
G}{\big(\mu^{-1}(x)\tensor e_{\mu} + e_{\mu}\tensor
\mu(x)\big)})=x$$ and
$$(\epsilon\tensor 1)\Delta(e_{\xi})= (\epsilon\tensor
1)(\sum_{\sigma\mu=\xi}{e_{\sigma}\tensor e_{\mu}} +
\Epsilon_{\xi})=
\sum_{\sigma\mu=\xi}{\delta_{1,\sigma}e_{\mu}}=e_{\xi}.$$ A
symmetric calculation shows that $(1 \otimes \epsilon)\Delta(x)=x$
and $(1 \otimes \epsilon)\Delta(e_{\xi})=e_{\xi}.$

The fact that $S$ is an antipode follows from the following
computations ($x\in L$):
$$\begin{array}{ll}
m(S\tensor 1)\Delta(x) & = m(S\tensor 1)({\displaystyle
\sum_{\mu\in
G}} {\big(\mu^{-1}(x) \tensor e_{\mu} + e_{\mu}\tensor \mu(x)\big)}) \vspace{3pt} \\
& = {\displaystyle \sum_{\mu\in G}}{\big(\mu^{-1}(x)e_{\mu} +
e_{\mu^{-1}}\mu(x)\big)} \vspace{3pt} \\
& =0 \vspace{3pt} \\
& =\epsilon(x)
\end{array}$$
and
$$\begin{array}{ll}
m(S\tensor 1)\Delta(e_{\xi}) & = m(S\tensor 1)({\displaystyle \sum_{\sigma\mu=\xi}}{e_{\sigma}\tensor e_{\mu}} + \Epsilon_{\xi}) \vspace{3pt} \\
& = {\displaystyle \sum_{\sigma\mu=\xi}}{e_{\sigma^{-1}}e_{\mu}} +
m(\Epsilon_{\xi}) \vspace{3pt} \\
& = \delta_{1,\xi}\big({\displaystyle \sum_{\sigma\in G}}{e_{\sigma}}+1_L \big) \vspace{3pt} \\
& = \epsilon(e_{\xi}).
\end{array}$$
The computations for $1\tensor S$ are exactly the same. This
completes the proof that $H$ is a Hopf algebra and the first part
of Theorem \ref{main3}.

We now prove that $H\otimes_k L\cong
\textrm{Fun}(\mathbb{Z}/2\mathbb{Z} \ltimes G,L)$. Since $L\tensor
_k L = \directsum_{\sigma\in G}L\Epsilon_{\sigma}$ and
$\textrm{Fun}(G,k)\otimes_k L \cong \textrm{Fun}(G,L)$, we know
that as an $L-$algebra $H\otimes_k L$ has a basis consisting of
$2n$ mutually orthogonal idempotents. The group
$Hom_{L-alg}(H\otimes_k L,L)$ thus contains $2n$ elements. There
is a natural mo\-no\-mor\-phism $Hom_{k-alg}(H,L) \rightarrow
Hom_{L-alg}(H\otimes_k L,L)$ given by extension of scalars. We
know that the group $Hom_{k-alg}(H,L)$ contains $2n$ elements and
is isomorphic to $\mathbb{Z}/2\mathbb{Z}\ltimes G$. The
mo\-no\-mor\-phism above is thus an isomorphism and therefore the
set of group-like elements of $(H\otimes_k L)^*$ is
$Hom_{L-alg}(H\otimes_k L,L)$, which is isomorphic to
$\mathbb{Z}_2\ltimes G$. Since $dim_L((H\otimes_k
L)^*)=dim_L(L[\mathbb{Z}/2\mathbb{Z}\ltimes G])$, it follows that
the Hopf algebras $H\otimes_k L$ and
$\textrm{Fun}(\mathbb{Z}_2\ltimes G,L)$ are isomorphic, as
claimed.
\end{proof}
\end{section}

\begin{section}{The proof of Theorem \ref{main1}}\label{sec-main-thm}
Let $L$ be a cyclic extension of $k$ with a cyclic Galois group
$G$, and let $\alpha \in Z^2(G,L^*)$. The crossed product algebra
$L^{\alpha}_tG$ has an $L$-basis given by
$\{U_{\sigma}\}_{\sigma\in G}$, and the multiplication is given by
the rule
$$(xU_{\sigma})(yU_{\tau})=x\sigma(y)\alpha(\sigma,\tau)U_{\sigma\tau},\textrm{
where }x,y\in L, \ \sigma,\tau\in G.$$ As mentioned in the
introduction, since the Galois group G is cyclic, we assume (as we
may) that $\alpha$ gets values in $k=L^G$ (in fact $\alpha$ is
cohomologous to a 2-cocycle $\beta$ of the form
$$\beta(\sigma^{i},\sigma^{j})=\left\{\begin{array}{ll}
1 & \textrm{if}\; i+j < |G| \\
b & \textrm{if}\; i+j\geq |G|
\end{array}\right., \textrm{ where } 0\leq i,j\leq |G|.$$
where $b$ is in $k=L^{G}$). In that case consider the
$k$-subalgebra of $L^{\alpha}_tG$ generated by the elements
$U_{\sigma}$. It is easy to see that this is the twisted group
algebra $k^{\alpha}G$ and by Theorem \ref{main2}, it is the
quotient of a semisimple Hopf algebra $A$ over $k$. By Theorem
\ref{main3}, we know that the field extension $L$ is a quotient of
a semisimple Hopf algebra $H$ over $k$. We now show how to
``amalgamate" these two constructions to obtain a Hopf algebra $X$
which projects onto the cyclic algebra $L^{\alpha}_tG$.

As a coalgebra, we let $X:=H\tensor A$. This means that the
comultiplication is given by the formula
\begin{equation}\label{define co-multiplication}\Delta(h\tensor
a)=\sum{h_{(1)}\tensor a_{(1)}\tensor h_{(2)}\tensor a_{(2)}}.
\end{equation}
The counit will be naturally given by
\begin{equation}\label{define counit} \epsilon(h\tensor
a)=\epsilon(h)\epsilon(a).
\end{equation}

Consider the algebra $L\tensor k^{\alpha^n}G$ with multiplication
\begin{equation}\label{define multiplication1}(x\tensor U^{(n)}_{\sigma})(y\tensor
U^{(n)}_{\tau})=x\sigma(y)\alpha^n(\sigma,\tau)\tensor
U^{(n)}_{\sigma\tau}.
\end{equation}
Notice that this is exactly the same as $L^{\alpha^n}_tG$.
Consider also the tensor algebra $Fun(G,k)\tensor k^{\alpha^n}G$,
whose multiplication is
\begin{equation}\label{define multiplication2}(e_{\mu}\tensor U^{(n)}_{\sigma})(e_{\nu}\tensor U^{(n)}_{\tau})=
\delta_{\mu,\nu}e_{\mu}\tensor
\alpha^n(\sigma,\tau)U^{(n)}_{\sigma\tau}.
\end{equation}
Then, as an algebra,
$$X=\left(\bigoplus_{n=0}^{m-1} L\tensor k^{\alpha^n}G \right) \bigoplus
\left(\bigoplus_{n=0}^{m-1} Fun(G,k) \tensor
k^{\alpha^n}G\right),$$ where $m$ is the order of $[\alpha]$.

The antipode $S$ on $X$ is given by
\begin{equation}\label{define antipode1}S(x\tensor U^{(n)}_{\sigma})=
\sigma^{-1}(x)\tensor
\alpha^n(\sigma,\sigma^{-1})/f(\sigma^{-1})U^{(m-n)}_{\sigma^{-1}},
\end{equation}
and
\begin{equation}\label{define antipode2}S(e_{\tau}\tensor U^{(n)}_{\sigma})
=e_{\tau^{-1}}\tensor
\alpha^n(\sigma,\sigma^{-1})/f(\sigma^{-1})U^{(m-n)}_{\sigma^{-1}}.
\end{equation}

\begin{proposition}
The above multiplication and antipode equip $X$ with a Hopf
algebra structure. The dimension of $X$ equals $2m|G|[L:k]$.
Furthermore, the crossed product $L^{\alpha}_tG$ is a quotient of
$X$.
\end{proposition}
\begin{proof}
Since the tensor product of two coalgebras is again a coalgebra,
we know that $X$ is a coalgebra. It is also easy to see that the
multiplication defined above makes $X$ into an algebra. So we only
need to check that the two structures are compatible. The fact
that $\epsilon$ is an algebra map is easy, and can be proved
directly, using Equations (\ref{define multiplication1}) and
(\ref{define multiplication2}). The fact that $\Delta$ is an
algebra map can also be proved directly by the above equations. We
give here a brief description of the computations. By the way
$\Delta$ was defined it is easy to see that
$\Delta(1_X)=\Delta(1\tensor1) = 1\tensor 1\tensor 1\tensor 1 =
1_{X\tensor X}$. Also, we have
$$\begin{array}{ll}
\Delta(h\tensor a) & = \sum{h_{(1)}\tensor
a_{(1)}\tensor h_{(2)}\tensor a_{(2)}} \vspace{3pt} \\
& = (\sum{h_{(1)}\tensor 1\tensor h_{(2)}\tensor
1})(\sum{1\tensor a_{(1)}\tensor 1\tensor a_{(2)}}) \vspace{3pt} \\
& = \Delta(h\tensor 1)\Delta(1\tensor a).
\end{array}$$
Since we already know that $\Delta$ is multiplicative on $H$ and
on $A$, it is easy to see that we only need to prove that for
every $h\in H$ and $a\in A$, $\Delta(1\tensor a)\Delta(h\tensor
1)= \Delta((1\tensor a)(h\tensor 1)).$ In order to show this we
need to recall the property about the idempotents of $L\tensor L$
previously used, that is,
$\omega_{\sigma}(\Epsilon_{\mu})=(\sigma\tensor\sigma)(\Epsilon_{\mu})=
\Epsilon_{\mu}$ for all $\sigma,\mu \in G.$ Writing
$\Epsilon_{\mu}=\sum_i{x_i\tensor y_i}$ this means that
$\sum_i{x_i\tensor y_i}=\sum_i{\sigma(x_i)\tensor \sigma(y_i)}$.

For $x\in L$ we compute:
\begin{eqnarray*}
\lefteqn{\Delta(1\tensor U^{(n)}_{\sigma})\Delta(x\tensor
1)}\\
& = & (\sum_{r+l\equiv n}{1\tensor \xi^{\sigma}_{r,l}
U^{(r)}_{\sigma}\tensor 1\tensor U^{(l)}_{\sigma}})(\sum_{\mu\in
G}{\mu^{-1}(x)\tensor 1 \tensor e_{\mu}\tensor 1 + e_{\mu}\tensor
1\tensor \mu(x)\tensor 1})\\
& = & \sum_{r+l\equiv n}\sum_{\mu\in
G}{\xi^{\sigma}_{r,l}\sigma\mu^{-1}(x)\tensor
U^{(r)}_{\sigma}\tensor e_{\mu}\tensor U^{(l)}_{\sigma} +
{\xi^{\sigma}_{r,l}}e_{\mu}\tensor U^{(r)}_{\sigma}\tensor
\sigma\mu(x)\tensor
U^{(l)}_{\sigma}} \\
& = & \sum_{r+l\equiv n}\sum_{\mu\in
G}{\xi^{\sigma}_{r,l}}\mu^{-1}\sigma(x)\tensor
U^{(r)}_{\sigma}\tensor e_{\mu}\tensor U^{(l)}_{\sigma}
+{\xi^{\sigma}_{r,l} e_{\mu}\tensor U^{(r)}_{\sigma}\tensor
\mu\sigma(x)\tensor
U^{(l)}_{\sigma}} \\
& = & \Delta(\sigma(x)\tensor U^{(n)}_{\sigma}) \\
& = & \Delta((1\tensor U^{(n)}_{\sigma})(x\tensor 1)).
\end{eqnarray*}
Let $\mu\in G$. Then
\begin{eqnarray*}
\lefteqn{ \Delta(1\tensor U^{(n)}_{\sigma})\Delta(e_{\mu}\tensor
1)}\\
& = & (\sum_{r+l\equiv n}{1\tensor \xi^{\sigma}_{r,l}
U^{(r)}_{\sigma}\tensor 1\tensor
U^{(l)}_{\sigma}})(\sum_{\rho\eta=\mu}{e_{\rho}\tensor 1\tensor
e_{\eta}\tensor 1} + \sum_i{x_i\tensor 1\tensor y_i\tensor 1})\\
& = & \sum_{r+l\equiv n}\sum_{\rho\eta=\mu}{e_{\rho}\tensor
\xi^{\sigma}_{r,l}U^{(r)}_{\sigma}\tensor e_{\eta}\tensor
U^{(l)}_{\sigma}} + \hspace{-3pt} \sum_{r+l\equiv
n}\sum_i{\sigma(x_i)\tensor\xi^{\sigma}_{r,l}U^{(r)}_{\sigma}\tensor\sigma(y_i)\tensor
U^{(l)}_{\sigma}}\\
& = & \sum_{r+l\equiv n}\sum_{\rho\eta=\mu}{e_{\rho}\tensor
\xi^{\sigma}_{r,l}U^{(r)}_{\sigma}\tensor e_{\eta}\tensor
U^{(l)}_{\sigma}} + \hspace{-3pt} \sum_{r+l\equiv
n}\sum_i{x_i\tensor\xi^{\sigma}_{r,l}U^{(r)}_{\sigma}\tensor
y_i\tensor U^{(l)}_{\sigma}}\\
& = & \Delta(e_{\mu}\tensor U^{(n)}_{\sigma}) \\
& = & \Delta((1\tensor U^{(n)}_{\sigma})(e_{\mu}\tensor 1)).
\end{eqnarray*}
This finishes the proof that $\Delta$ is multiplicative.

Finally we show that $S$ is an antipode. Notice that by the way
$S$ was defined, we have $S(h\tensor a)= (1\tensor
S_A(a))(S_H(h)\tensor 1)$. This can be easily checked using
equations (\ref{define multiplication1}) -- (\ref{define
antipode2}). We now compute
$$\begin{array}{ll}
m(S\tensor 1)\Delta(h\tensor a) & = m(S\tensor
1)({\displaystyle \sum}{h_{(1)}\tensor a_{(1)}\tensor h_{(2)}\tensor a_{(2)}}) \vspace{3pt} \\
& = {\displaystyle \sum}{(1\tensor
S_A(a_{(1)}))(S_H(h_{(1)})\tensor
1)h_{(2)}\tensor a_{(2)}} \vspace{3pt} \\
& =\epsilon(h){\displaystyle \sum}{(1\tensor S_A(a_{(1)}))(1\tensor a_{(2)})} \vspace{3pt} \\
& = \epsilon(h)\epsilon(a)1\tensor 1\tensor 1\tensor 1\\
& = \epsilon(h\tensor a)1_X.
\end{array}$$
Similarly, $m(1 \tensor S)\Delta(h\tensor a)=\epsilon(h\tensor
a)1_X.$ Thus $S$ is an antipode for $X$, and $X$ is indeed a Hopf
algebra. Since $L^\alpha_tG$ is a direct summand of $X$, it is a
quotient of $X$. The proof is complete.
\end{proof}

\begin{remark}
One can also view $X$ as a bicrossed product of the Hopf algebras
$A$ and $H$. By Theorems \ref{main2} and \ref{main3}, $X$ is a
$\bar{k}-$form of the bicrossed product of the Hopf algebras
$\bar{k}\widehat{G}$ and $\text{Fun}(\mathbb{Z}/2\mathbb{Z}
\ltimes G,\bar{k})$. (See e.g. \cite[IX.2]{K} for details about
bicrossed products.)
\end{remark}

\begin{remark} In order to construct $X$ the field extension $L$ need not by cyclic.
All we really need is that $L$ will be an abelian extension, and
that the cocycle $\alpha$ will be in the image of
$H^2(G,k^*)\rightarrow H^2(G,L^*)$.
\end{remark}
\end{section}
\vspace{0.5cm}

\section*{Acknowledgments}
The second named author's research was supported by projects
MTM2005-03227 from MEC and FEDER and P06-FQM-01889 from Junta de
Andaluc\'{\i}a.
%
The third named author's research was supported by the Israel
Science Foundation (grant No. 125/05).

\vspace{0.5cm}


\begin{thebibliography}{ABCD}

\bibitem[AS1]{AS1}
E. Aljadeff and J. Sonn, Projective Schur algebras have abelian
splitting fields. J. Algebra 175 (1995), no. 1, 179--187.

\bibitem[AS2]{AS2}
E. Aljadeff and J. Sonn, On the projective Schur group of a field.
J. Algebra 178 (1995), no. 2, 530--540.

\bibitem[FS]{FS}
B. Fein and M. Schacher, Brauer groups of rational function
fields, in "Groupe de Brauer", Lecture Notes in Mathematics, vol
844, Springer-Verlag, New York/Berlin, 1981.

\bibitem[G]{G}
S. Gelaki, Semisimple triangular Hopf algebras and Tannakian
categories, {\em Proceedings of Symposia in Pure Mathematics} {\bf
70} (2002), 497--516 (editors M. Fried and Y. Ihara, 1999 Von
Neumann Conference on Arithmetic Fundamental Groups and
Noncommutative Algebra, August 16-27, 1999 MSRI).

\bibitem[J]{J}
G. Janusz, The Schur group of an algebraic number field, {\em Ann.
of Math.} (2) {\bf 103} (1976), no. 2, 253--281.

\bibitem[K]{K}
C. Kassel, Quantum Groups, Graduate Texts in Mathematics {\bf
155}, Springer-Verlag, 1995.

\bibitem[Y]{Y} T. Yamada, The Schur Subgroup of the Brauer Group,
Lecture Notes in Mathematics {\bf 397}, Springer-Verlag, 1970.
\end{thebibliography}

\end{document}